\documentclass[12pt]{article}
\usepackage{amsfonts}
\usepackage{a4wide}
\usepackage{amsmath,amscd,amsthm,a4,amssymb}
\usepackage{sectsty}

\sectionfont{\centering}
\numberwithin{equation}{section}

\begin{document}

\begin{center}
{\large \textbf{A REPRESENTATION PROBLEM FOR SMOOTH SUMS OF RIDGE FUNCTIONS}}

\

\textbf{Rashid A. Aliev} \

Institute of Mathematics and Mechanics, NAS of Azerbaijan, Baku, Azerbaijan

Baku State University, Baku, Azerbaijan

{e-mail:} {aliyevrashid@mail.ru}

\ \

\textbf{Vugar E. Ismailov}

Institute of Mathematics and Mechanics, NAS of Azerbaijan, Baku, Azerbaijan

{e-mail:} {vugaris@mail.ru}
\end{center}

\bigskip

\textbf{Abstract.} In this paper we prove that if a multivariate function
of a certain smoothness class is represented by a sum of $k$ arbitrarily
behaved ridge functions, then it can be represented by a sum of $k$ ridge
functions of the same smoothness class and a polynomial of degree at most $k-1$.
This solves the problem posed by A. Pinkus in his monograph ``Ridge
Functions" up to a multivariate polynomial.

\bigskip

\textbf{Mathematics Subject Classification:} 26B40, 39B22.

\textbf{Keywords:} ridge function; Cauchy functional equation; difference
property; polynomial function.

\bigskip

\bigskip

\begin{center}
{\large \textbf{1. Introduction}}
\end{center}

This paper explores the following open question raised in Buhmann and Pinkus
\cite{12}, and Pinkus \cite[p. 14]{Pinkus1}. Assume we are given a function $%
f(\mathbf{x})=f(x_{1},...,x_{n})$ of the form
\begin{equation*}
f(\mathbf{x})=\sum_{i=1}^{k}f_{i}(\mathbf{a}^{i}\cdot \mathbf{x}),\eqno(1.1)
\end{equation*}%
where the $\mathbf{a}^{i},$ $i=1,...,k,$ are pairwise linearly independent
vectors (directions) in $\mathbb{R}^{n}$, $f_{i}$ are arbitrarily behaved
univariate functions and $\mathbf{a}^{i}\cdot \mathbf{x}$ are standard inner
products. Assume, in addition, that $f$ is of a certain smoothness class,
that is, $f\in C^{s}(\mathbb{R}^{n})$, where $s\geq 0$ (with the convention
that $C^{0}(\mathbb{R}^{n})=C(\mathbb{R}^{n})$). Is it true that there will
always exist $g_{i}\in C^{s}(\mathbb{R})$ such that
\begin{equation*}
f(\mathbf{x})=\sum_{i=1}^{k}g_{i}(\mathbf{a}^{i}\cdot \mathbf{x})\text{ ?}%
\eqno(1.2)
\end{equation*}

Functions of the form $g(\mathbf{a}\cdot \mathbf{x})$, involved in the right
hand sides of (1.1) and (1.2), are called ridge functions. These functions
appear in various fields and under various guises. They appear in partial
differential equations (where they are called \textit{plane waves}, see,
e.g., \cite{69}), in computerized tomography (see, e.g., \cite{97,111}), in
statistics (especially, in the theory of projection pursuit and projection
regression; see, e.g., \cite{27,33}). Ridge functions are also the
underpinnings of many central models in neural networks which has become
increasingly more popular in the last few decades in many fields of science
and engineering (see \cite{119} and a great deal of references therein).
Finally, these functions are used in modern approximation theory as an
effective and convenient tool for approximating complicated multivariate
functions (see, e.g., \cite{49,50,54,88,95,102}). We refer the reader to the
monograph by Pinkus \cite{Pinkus1} for a detailed and systematic study of
ridge functions.

Note that for $k=1$ and $k=2$ the above problem is easily solved. Indeed for
$k=1$ by choosing $\mathbf{c}\in \mathbb{R}^{n}$ satisfying $\mathbf{a}%
^{1}\cdot \mathbf{c}=1$, we have that $f_{1}(t)=f(t\mathbf{c)}$ is in $C^{s}(%
\mathbb{R})$. The same argument can be carried out for the case $k=2.$ In
this case, since the vectors $\mathbf{a}^{1}$ and $\mathbf{a}^{2}$ are
linearly independent, there exists a vector $\mathbf{c}\in \mathbb{R}^{n}$
satisfying $\mathbf{a}^{1}\cdot \mathbf{c}=1$ and $\mathbf{a}^{2}\cdot
\mathbf{c}=0.$ Therefore, we obtain that the function $f_{1}(t)=f(t\mathbf{c)%
}-f_{2}(0)$ is in the class $C^{s}(\mathbb{R})$. Similarly, one can verify
that $f_{2}\in C^{s}(\mathbb{R})$.

The above cases with one and two ridge functions in (1.1) show that the
functions $f_{i}$ inherit smoothness properties of the given $f$. The
picture is absolutely different if the number of directions $k\geq 3$. For $%
k=3$, there are ultimately smooth functions which decompose into sums of
very badly behaved ridge functions. This phenomena comes from the classical
Cauchy Functional Equation (CFE). This equation,%
\begin{equation*}
h(x+y)=h(x)+h(y),\text{ }h:\mathbb{R\rightarrow R}\text{,}\eqno(1.3)
\end{equation*}%
looks very simple and has a class of simple solutions $h(x)=cx,$ $c\in
\mathbb{R}$. However, it easily follows from Hamel basis theory that CFE
also has a large class of wild solutions. These solutions are called
\textquotedblleft wild" because they are extremely pathological. They are,
for example, not continuous at a point, not monotone on an interval, not
bounded on any set of positive measure (see, e.g., \cite{1}). Let $h_{1}$ be
any wild solution of the equation (1.3). Then the zero function can be
represented as%
\begin{equation*}
0=h_{1}(x)+h_{1}(y)-h_{1}(x+y).\eqno(1.4)
\end{equation*}%
Note that the functions involved in (1.4) are bivariate ridge functions with
the directions $(1,0)$, $(0,1)$ and $(1,1)$, respectively. This example
shows that for $k\geq 3$ the functions $f_{i}$ in (1.1) may not inherit
smoothness properties of the function $f$, which in the case of (1.4) is the
identically zero function. Thus the above problem arises naturally.

However, it was shown by some authors that, additional conditions on $f_{i}$
or the directions $\mathbf{a}^{i}$ guarantee smoothness of the
representation (1.1). It was first proved by Buhmann and Pinkus \cite{12}
that if in (1.1) $f\in C^{s}(\mathbb{R}^{n})$, $s\geq k-1$ and $f_{i}\in
L_{loc}^{1}(\mathbb{R)}$ for each $i$, then $f_{i}\in C^{s}(\mathbb{R)}$ for
$i=1,...,k.$ Later Pinkus \cite{Pinkus} found a strong relationship between
CFE and the problem of smoothness in ridge function representation. He
generalized extensively the previous result of Buhmann and Pinkus \cite{12}.
He showed that the solution is quite simple and natural if the functions $%
f_{i}$ are taken from a certain class $\mathcal{B}$ of real-valued functions
defined on $\mathbb{R}$. $\mathcal{B}$ includes, for example, the set of
continuous functions, the set of bounded functions, the set of Lebesgue
measurable functions (for the precise definition of $\mathcal{B}$ see
Section 4). The result of Pinkus states that if in (1.1) $f\in C^{s}(\mathbb{%
R}^{n})$ and each $f_{i}\in \mathcal{B}$, then necessarily $f_{i}\in C^{s}(%
\mathbb{R)}$ for $i=1,...,k$ (see \cite{Pinkus}).

Note that severe restrictions on the directions $\mathbf{a}^{i}$ also
guarantee smoothness of the representation (1.1). For example, in \cite{K1}
it was easily proven that in (1.1) the inclusions $f_{i}\in C^{s}(\mathbb{R}%
) $, $i=1,...,k,$ are automatically valid if the directions $\mathbf{a}^{i}$
are linearly independent and if these directions are not linearly
independent, then there exists $f\in C^{s}(\mathbb{R}^{n})$ of the form
(1.1) such that the $f_{i}\notin C^{s}(\mathbb{R}),$ $i=1,...,k.$

The above result of Pinkus was a starting point for further research on
continuous and smooth sums of ridge functions. Much work in this direction
was done by Konyagin and Kuleshov \cite{K1,K2}, and Kuleshov \cite{K4}. They
mainly analyze the continuity of $f_{i}$, that is, the question of if and
when continuity of $f$ guarantees the continuity of $f_{i}$. There are also
other results concerning different properties, rather than continuity, of $%
f_{i}$. Most results in \cite{K1,K2,K4} involve certain subsets (convex open
sets, convex bodies, etc.) of $\mathbb{R}^{n}$ instead of only $\mathbb{R}%
^{n}$ itself.

In \cite{A1}, we gave a partial solution to the above representation
problem. Our solution comprises the cases in which $s\geq 1$ and $k-1$
directions of the given $k$ directions are linearly independent. For
bivariate functions having degree of smoothness $s\geq k-2,$ the problem was
solved in \cite{A2}.

Kuleshov \cite{K3} generalized our result \cite[Theorem 2.3]{A1} to all
possible cases of $s$. That is, he proved that if a function $f\in C^{s}(%
\mathbb{R}^{n})$, where $s\geq 0$, is of the form (1.1) and $(k-1)$-tuple of
the given set of $k$ directions $\mathbf{a}^{i}$ forms a linearly
independent system, then there exist $g_{i}\in C^{s}(\mathbb{R})$, $%
i=1,...,k $, such that (1.2) holds (see \cite[Theorem 3]{K3}). In \cite{A0},
we reproved this result using completely different ideas. Note that our
proof contains a theoretical method for constructing the functions $g_{i}\in
C^{s}(\mathbb{R})$ in (1.2) (see \cite[Theorem 2.1, Theorem 2.2]{A0}). Using
this method, we also estimated the modulus of continuity of $f_{i}$ in terms
of the modulus of continuity of $f$ (see \cite[Remark 2]{A0}).

In this paper, based on the theory of polynomial functions (see \cite[%
Section 15.9]{Kuc}), we give a solution to the above representation problem
up to some multivariate polynomial. That is, we show that if (1.1) holds for
$f\in C^{s}(\mathbb{R}^{n})$ and arbitrarily behaved $f_{i}$, then there
exist $g_{i}\in C^{s}(\mathbb{R})$ such that
\begin{equation*}
f(\mathbf{x})=\sum_{i=1}^{k}g_{i}(\mathbf{a}^{i}\cdot \mathbf{x})+P(\mathbf{x%
}),
\end{equation*}%
where $P(\mathbf{x})$ is a polynomial of degree at most $k-1$. This leads to
a complete solution of the problem in the case when the space dimension $n=2$.
We also prove that if the directions $\mathbf{a}^{i}$ in (1.1) have only
rational coordinates, then the above polynomial term does not appear, and
hence (1.2) holds.

\bigskip

\bigskip

\bigskip

\bigskip

\bigskip

\begin{center}
{\large \textbf{2. Polynomial functions of k-th order}}
\end{center}

Given $h_{1},...,h_{k}\in \mathbb{R}$, we define inductively the difference
operator $\Delta _{h_{1}...h_{k}}$ as follows
\begin{eqnarray*}
\Delta _{h_{1}}f(x) &:&=f(x+h_{1})-f(x), \\
\Delta _{h_{1}...h_{k}}f &:&=\Delta _{h_{k}}(\Delta _{h_{1}...h_{k-1}}f),%
\text{ }f:\mathbb{R\rightarrow R}.
\end{eqnarray*}%
If $h_{1}=\cdots=h_{k}=h,$ then we write briefly $\Delta _{h}^{k}f$ instead
of $\Delta _{\underset{n\text{ times}}{\underbrace{h...h}}}f$. For various
properties of difference operators see \cite[Section 15.1]{Kuc}.

\bigskip

\textbf{Definition 2.1 }(see \cite{Kuc}). A function $f:\mathbb{R\rightarrow
R}$ is called a \textit{polynomial function} of order $k$ ($k\in \mathbb{N}$%
) if for every $x\in \mathbb{R}$ and $h\in \mathbb{R}$ we have
\begin{equation*}
\Delta _{h}^{k+1}f(x)=0.
\end{equation*}

It can be shown that if $\Delta _{h}^{k+1}f=0$ for any $h\in \mathbb{R}$,
then $\Delta _{h_{1}...h_{k+1}}f=0$ for any $h_{1},...,h_{k+1}\in \mathbb{R}$
(see \cite[Theorem 15.3.3]{Kuc}). A polynomial of degree at most $k$ is a
polynomial function of order $k$ (see \cite[Theorem 15.9.4]{Kuc}). The
polynomial functions generalize ordinary polynomials, and reduce to the
latter under mild regularity assumptions. For example, if a polynomial
function is continuous at one point, or bounded on a set of positive
measure, then it continuous at all points (see \cite{Cies, Kurepa}), and
therefore is a polynomial of degree $k$ (see \cite[Theorem 15.9.4]{Kuc}).

Basic results concerning polynomial functions are due to S. Mazur-W. Orlicz
\cite{Maz}, McKiernan \cite{Mc}, Djokovi\'{c} \cite{Djok}. The following
theorem, which we will use in the sequel, yield implicitly the general
construction of polynomial functions.

\bigskip

\textbf{Theorem 2.1 }(see \cite[Theorems 15.9.1 and 15.9.2]{Kuc}). \textit{A
function $f:\mathbb{R\rightarrow R}$ is a polynomial function of order $k$
if and only if it admits a representation}

\begin{equation*}
f=f_{0}+f_{1}+...+f_{k},
\end{equation*}%
\textit{where $f_{0}$ is a constant and $f_{j}:\mathbb{R\rightarrow R}$, $%
j=1,...,k$, are diagonalizations of $j$-additive symmetric functions $F_{j}:%
\mathbb{R}^{j}\mathbb{\rightarrow R}$, i.e.,}

\begin{equation*}
f_{j}(x)=F_{j}(x,...,x).
\end{equation*}

\bigskip

Note that a function $F_{p}:\mathbb{R}^{p}\mathbb{\rightarrow R}$ is called $%
p$-additive if for every $j,$ $1\leq j\leq p,$ and for every $%
x_{1},...,x_{p},y_{j}\in \mathbb{R}$

\begin{equation*}
F(x_{1},...,x_{j}+y_{j},...,x_{p})=F(x_{1},...,x_{p})+F(x_{1},...,x_{j-1},y_{j},x_{j+1},...,x_{p}),
\end{equation*}%
i.e., $F$ is additive in each of its variables $x_{j}$ (see \cite[p. 363]%
{Kuc}). A simple example of a $p$-additive function is given by the product

\begin{equation*}
f_{1}(x_{1})\times\cdots\times f_{p}(x_{p}),
\end{equation*}%
where the univariate functions $f_{j},$ $j=1,...,p$, are additive.

Following de Bruijn, we say that a class $\mathcal{D}$ of real functions has
the \textit{difference property} if any function $f:\mathbb{R\rightarrow R}$
such that $\bigtriangleup _{h}f\in \mathcal{D}$ for all $h\in \mathbb{R}$,
admits a decomposition $f=g+S$, where $g\in \mathcal{D}$ and $S$ satisfies
the Cauchy Functional Equation (1.3). Several classes with the difference
property are investigated in de Bruijn \cite{B1,B2}. Some of these classes
are:

1) $C(\mathbb{R)}$, continuous functions;

2) $C^{s}(\mathbb{R)}$, functions with continuous derivatives up to order $s$%
;

3) $C^{\infty }(\mathbb{R)}$, infinitely differentiable functions;

4) analytic functions;

5) functions which are absolutely continuous on any finite interval;

6) functions having bounded variation over any finite interval;

7) algebraic polynomials;

8) trigonometric polynomials;

9) Riemann integrable functions.

\bigskip

A natural generalization of classes with the difference property are classes
of functions with the difference property of $k$-th order.

\bigskip

\textbf{Definition 2.2 }(see \cite{Gajda}). A class $\mathcal{F}$ is said to
have the \textit{difference property of }$k$\textit{-th order} if any
function $f:\mathbb{R\rightarrow R}$ such that $\bigtriangleup _{h}^{k}f\in
\mathcal{F}$ for all $h\in \mathbb{R}$, admits a decomposition $f=g+H$,
where $g\in \mathcal{F}$ and $H$ is a polynomial function of $k$-th order.

\bigskip

It is not difficult to see that the class $\mathcal{F}$ has the difference
property of first order if and only if it has the difference property in de
Bruijn's sense. There arises a natural question: which of the above classes
have difference properties of higher orders? Gajda \cite{Gajda} considered
this question in its general form, for functions defined on a locally
compact Abelian group and showed that for any $k\in \mathbb{N}$, continuous
functions have the difference property of $k$-th order (see \cite[Theorem 4]%
{Gajda}). The proof of this result is based on several lemmas, in
particular, on the following lemma, which we will also use in the sequel.

\bigskip

\textbf{Lemma 2.1.} (see \cite[Lemma 5]{Gajda}). \textit{For each $k\in
\mathbb{N}$ the class of all continuous functions defined on $\mathbb{R}$
has the difference property of $k$-th order.}

\bigskip

In fact, Gajda \cite{Gajda} proved this lemma for Banach space valued
functions, but the simplest case with the space $\mathbb{R}$ has all the
difficulties. Unfortunately, the proof of the lemma has an essential gap.
The author of \cite{Gajda} tried to reduce the proof to $\mod1$ periodic
functions, but made a mistake in proving the continuity of the difference $%
\Delta _{h_{1}...h_{k-1}}(f-f^{\ast })$. Here $f^{\ast }:\mathbb{%
R\rightarrow R}$ is a $\mod1$ periodic function defined on the interval $%
[0,1)$ as $f^{\ast }(x)=f(x)$ and extended to the whole $\mathbb{R}$ with
the period $1$. That is, $f^{\ast }(x)=f(x)$ for $x\in \lbrack 0,1)$ and $%
f^{\ast }(x+1)=f^{\ast }(x)$ for $x\in \mathbb{R}$. In the proof, the author
of \cite{Gajda} takes a point $x\in \lbrack m,m+1)$ and writes that

\begin{equation*}
\Delta _{h_{1}...h_{k-1}}(f-f^{\ast })(x)=\Delta
_{h_{1}...h_{k-1}}(f(x)-f(x-m))\text{,}
\end{equation*}%
which is not valid. Even though $f^{\ast }(x)=f(x-m)$ for any $x\in \lbrack
m,m+1)$, the differences $\Delta _{h_{1}...h_{k-1}}f^{\ast }(x)$ and $\Delta
_{h_{1}...h_{k-1}}f(x-m)$ are completely different, since the latter may
involve values of $f$ at points outside $[0,1)$, which have no relationship
with the definition of $f^{\ast }$.

In the next section, we give a new proof for Lemma 2.1 (see Theorem 3.1). We
hope that our proof is free from mathematical errors and thus the above
lemma itself is valid.

\bigskip

\bigskip

\bigskip

\bigskip

\begin{center}
{\large \textbf{3. Some auxiliary results on polynomial functions}}
\end{center}

In this section, we do further research on polynomial functions and prove
some auxiliary results.

\bigskip

\textbf{Lemma 3.1.} \textit{If $f:\mathbb{R\rightarrow R}$ is a polynomial
function of order $k$, then for any $p\in $ $\mathbb{N}$ and any fixed $\xi
_{1},...,\xi _{p}\in \mathbb{R}$, the function}%
\begin{equation*}
g(x_{1},...,x_{p})=f(\xi _{1}x_{1}+\cdots+\xi _{p}x_{p}),
\end{equation*}
\textit{considered on the $p$ dimensional space $\mathbb{Q}^{p}$ of rational
vectors, is an ordinary polynomial of degree at most $k$.}

\bigskip

\textbf{Proof.} By Theorem 2.1,

\begin{equation*}
f=\sum_{m=0}^{k}f_{m},\eqno(3.1)
\end{equation*}%
where $f_{0}$ is a constant and $f_{m}:\mathbb{R\rightarrow R}$, $1,...,m$,
are diagonalizations of $m$-additive symmetric functions $F_{m}:\mathbb{R}%
^{m}\mathbb{\rightarrow R}$, i.e.,

\begin{equation*}
f_{m}(x)=F_{m}(x,...,x).
\end{equation*}%
For a $m$-additive function $F_{m}$ the equality

\begin{equation*}
F_{m}(\xi _{1},...,\xi _{i-1},r\xi _{i},\xi _{i+1},...,\xi _{m})=rF_{m}(\xi
_{1},...,\xi _{m})
\end{equation*}%
holds for all $i=1,...,m$ and any $r\in \mathbb{Q}$, $\xi _{i}\in $ $\mathbb{%
R}$, $i=1,...,m$ (see \cite[Theorem 13.4.1]{Kuc}). Using this, it is not
difficult to verify that for any $(x_{1},...,x_{p})\in \mathbb{Q}^{p}$,

\begin{eqnarray*}
f_{m}(\xi _{1}x_{1}+\cdots+\xi _{p}x_{p}) &=&F_{m}(\xi _{1}x_{1}+\cdots+\xi
_{p}x_{p},...,\xi _{1}x_{1}+\cdots+\xi _{p}x_{p}) \\
&=&\sum_{\substack{ 0\leq s_{i}\leq m,~\overline{i=1,p}  \\ %
s_{1}+\cdots+s_{p}=m }}A_{s_{1}...s_{p}}F_{m}(\underset{s_{1}}{\underbrace{%
\xi _{1},...,\xi _{1}}},...,\underset{s_{p}}{\underbrace{\xi _{p},...,\xi
_{p}}})x_{1}^{s_{1}}...x_{p}^{s_{p}}.
\end{eqnarray*}%
Here $A_{s_{1}...s_{p}}$ are some coefficients, namely $%
A_{s_{1}...s_{p}}=m!/(s_{1}!...s_{p}!).$ Considering the last formula in
(3.1), we conclude that the function $g(x_{1},...,x_{p})$, restricted to $%
\mathbb{Q}^{p}$, is a polynomial of degree at most $k$.

\bigskip

\textbf{Lemma 3.2.} \textit{Assume $f$ is a polynomial function of order $k$%
. Then there exists a polynomial function $H$ of order $k+1$ such that $%
H(0)=0$ and}

\begin{equation*}
f(x)=H(x+1)-H(x).\eqno(3.2)
\end{equation*}

\bigskip

\textbf{Proof.} Consider the function

\begin{equation*}
H(x):=xf(x)+\sum_{i=1}^{k}(-1)^{i}\frac{x(x+1)...(x+i)}{(i+1)!}\Delta
_{1}^{i}f(x).\eqno(3.3)
\end{equation*}%
Clearly, $H(0)=0.$ We are going to prove that $H$ is a polynomial function
of order $k+1$ and satisfies (3.2).

Let us first show that for any polynomial function $g$ of order $m$ the
function $G_{1}(x)=xg(x)$ is a polynomial function of order $m+1.$ Indeed,
for any $h_{1},...,h_{m+2}\in \mathbb{R}$ we can write that

\begin{equation*}
\Delta _{h_{1}...h_{m+2}}G_{1}(x)=(x+h_{1}+\cdots+h_{m+2})\Delta
_{h_{1}...h_{m+2}}g(x)+\sum_{i=1}^{m+2}h_{i}\Delta
_{h_{1}...h_{i-1}h_{i+1...}h_{m+2}}g(x).\eqno(3.4)
\end{equation*}%
The last formula is verified directly by using the known product property of
differences, that is, the equality

\begin{equation*}
\Delta _{h}(g_{1}g_{2})=g_{1}\Delta _{h}g_{2}+g_{2}\Delta _{h}g_{1}+\Delta
_{h}g_{1}\Delta _{h}g_{2}.\eqno(3.5)
\end{equation*}%
Now since $g$ is a polynomial function of order $m$, all summands in (3.4)
is equal to zero; hence we obtain that $G_{1}(x)$ is a polynomial function
of order $m+1$. By induction, we can prove that the function $%
G_{p}(x)=x^{p}g(x)$ is a polynomial function of order $m+p.$ Since $\Delta
_{1}^{i}f(x)$ in (3.3) is a polynomial function of order $k-i$, it follows
that all summands in (3.3) are polynomial functions of order $k+1$.
Therefore, $H(x)$ is a polynomial function of order $k+1$.

Now let us prove (3.2). Considering the property (3.5) in (3.3) we can write
that%
\begin{equation*}
\Delta _{1}H(x)=\left[ f(x)+(x+1)\Delta _{1}f(x)\right]
\end{equation*}
\begin{equation*}
+\sum_{i=1}^{k}(-1)^{i}\left[ \frac{(x+1)...(x+i+1)}{(i+1)!}\Delta
_{1}^{i+1}f(x)+\Delta _{1}\left( \frac{x(x+1)...(x+i)}{(i+1)!}\right) \Delta
_{1}^{i}f(x)\right] .\eqno(3.6)
\end{equation*}

Note that in (3.6)

\begin{equation*}
\Delta _{1}\left( \frac{x(x+1)...(x+i)}{(i+1)!}\right) =\frac{(x+1)...(x+i)}{%
i!}.
\end{equation*}%
Considering this and the assumption $\Delta _{1}^{k+1}f(x)=0$, it follows
from (3.6) that

\begin{equation*}
\Delta _{1}H(x)=f(x),
\end{equation*}%
that is, (3.2) holds.

\bigskip

The next lemma is due to Gajda \cite{Gajda}.

\bigskip

\textbf{Lemma 3.3 }(see \cite[Corollary 1]{Gajda}). \textit{Let $f:$ $%
\mathbb{R\rightarrow R}$ be a $\mod1$ periodic function such that, for any $%
h_{1},...,h_{k}\in \mathbb{R}$, $\Delta _{h_{1}...h_{k}}f$ is continuous.
Then there exist a continuous function $g:$ $\mathbb{R\rightarrow R}$ and a
polynomial function $H$ of $k$-th order such that $f=g+H$.}

\bigskip

The following theorem generalizes de Bruijn's theorem (see \cite[Theorem 1.1]%
{B1}) on the difference property of continuous functions and shows that
Gajda's above lemma (see Lemma 2.1) is valid. Note that the main result of
\cite{Gajda} also uses this theorem.

\bigskip

\textbf{Theorem 3.1.} \textit{Assume for any $h_{1},...,h_{k}\in \mathbb{R}$%
, the difference $\Delta _{h_{1}...h_{k}}f(x)$ is a continuous function of
the variable $x$. Then there exist a function $g\in C(\mathbb{R})$ and a
polynomial function $H$ of $k$-th order with the property $H(0)=0$ such that}

\begin{equation*}
f=g+H.
\end{equation*}

\bigskip

\textbf{Proof.} We prove this theorem by induction. For $k=1$, the theorem
is the result of de Bruijn: if $f$ is such that, for each $h$, $\Delta
_{h}f(x)$ is a continuous function of $x$, then it can be written in the
form $g+H$, where $g$ is continuous and $H$ is additive (that is, satisfies
the Cauchy Functional Equation). Assume that the theorem is valid for $k-1.$
Let us prove it for $k$. Without loss of generality we may assume that $%
f(0)=f(1)$. Otherwise, we can prove the theorem for $f_{0}(x)=f(x)-\left[
f(1)-f(0)\right] x$ and then automatically obtain its validity for $f$.

Consider the function

\begin{equation*}
F_{1}(x)=f(x+1)-f(x)\text{, }x\in \mathbb{R}.\eqno(3.7)
\end{equation*}%
Since for\ any $h_{1},...,h_{k}\in \mathbb{R}$, $\Delta _{h_{1}...h_{k}}f(x)$
is a continuous function of $x$ and $\Delta _{h_{1}...h_{k-1}}F_{1}=\Delta
_{h_{1}...h_{k-1}1}f$, the difference $\Delta _{h_{1}...h_{k-1}}F_{1}(x)$
will be a continuous function of $x$, as well. By assumption, there exist a
function $g_{1}\in C(\mathbb{R})$ and a polynomial function $H_{1}$ of $%
(k-1) $-th order with the property $H_{1}(0)=0$ such that

\begin{equation*}
F_{1}=g_{1}+H_{1}.\eqno(3.8)
\end{equation*}%
It follows from Lemma 3.2 that there exists a polynomial function $H_{2}$ of
order $k$ such that $H_{2}(0)=0$ and

\begin{equation*}
H_{1}(x)=H_{2}(x+1)-H_{2}(x).\eqno(3.9)
\end{equation*}%
Substituting (3.9) in (3.8) we obtain that

\begin{equation*}
F_{1}(x)=g_{1}(x)+H_{2}(x+1)-H_{2}(x).\eqno(3.10)
\end{equation*}%
It follows from (3.7) and (3.10) that

\begin{equation*}
g_{1}(x)=\left[ f(x+1)-H_{2}(x+1)\right] -\left[ f(x)-H_{2}(x)\right] .\eqno%
(3.11)
\end{equation*}

Consider the function

\begin{equation*}
F_{2}=f-H_{2}.\eqno(3.12)
\end{equation*}%
Since $H_{2}$ is a polynomial function of order $k$ and for any $%
h_{1},...,h_{k}\in \mathbb{R}$ the difference $\Delta _{h_{1}...h_{k}}f(x)$
is a continuous function of $x$, we obtain that $\Delta
_{h_{1}...h_{k}}F_{2}(x)$ is also a continuous function of $x$. In addition,
since $f(0)=f(1)$ and $H_{2}(0)=H_{2}(1)=0$, it follows from (3.12) that $%
F_{2}(0)=F_{2}(1)$. We will use these properties of $F_{2}$ below.

Let us write (3.11) in the form

\begin{equation*}
g_{1}(x)=F_{2}(x+1)-F_{2}(x),\eqno(3.13)
\end{equation*}%
and define the following $\mod 1$ periodic function

\begin{eqnarray*}
F^{\ast }(x) &=&F_{2}(x)\text{ for }x\in \lbrack 0,1), \\
F^{\ast }(x+1) &=&F^{\ast }(x)\text{ for }x\in \mathbb{R}.
\end{eqnarray*}

Consider the function

\begin{equation*}
F=F_{2}-F^{\ast }.\eqno(3.14)
\end{equation*}%
Let us show that $F\in C(\mathbb{R)}$. Indeed since $F(x)=0$ for $x\in
\lbrack 0,1)$, $F$ is continuous on $(0,1)$. Consider now the interval $%
[1,2) $. For any $x\in \lbrack 1,2)$ by the definition of $F^{\ast }$ and
(3.13) we can write that

\begin{equation*}
F(x)=F_{2}(x)-F_{2}(x-1)=g_{1}(x-1).\eqno(3.15)
\end{equation*}%
Since $g_{1}\in C(\mathbb{R)}$, it follows from (3.15) that $F$ is
continuous on $(1,2)$. Note that by (3.13) $g_{1}(0)=0$; hence $%
F(1)=g_{1}(0)=0$. Since $F\equiv 0$ on $[0,1)$, $F(1)=0$ and $F\in C(1,2),$
we obtain that $F$ is continuous on $(0,2)$. Consider the interval $[2,3)$.
For any $x\in \lbrack 2,3)$ we can write that

\begin{equation*}
F(x)=F_{2}(x)-F_{2}(x-2)=g_{1}(x-1)+g_{1}(x-2).\eqno(3.16)
\end{equation*}%
Since $g_{1}\in C(\mathbb{R)}$, $F$ is continuous on $(2,3)$. Note that by
(3.15) $\lim_{x\rightarrow 2-}F(x)=g_{1}(1)$ and by (3.16) $F(2)=g_{1}(1).$
We obtain from these arguments that $F$ is continuous on $(0,3)$. In the
same way, we can prove that $F$ is continuous on $(0,m)$ for any $m\in
\mathbb{N}$.

Similar arguments can be used to prove the continuity of $F$ on $(-m,0)$ for
any $m\in \mathbb{N}$. We show it for the first interval $[-1,0)$. For any $%
x\in \lbrack -1,0)$ by the definition of $F^{\ast }$ and (3.13) we can write
that

\begin{equation*}
F(x)=F_{2}(x)-F_{2}(x+1)=-g_{1}(x).
\end{equation*}%
Since $g_{1}\in C(\mathbb{R)}$, it follows that $F$ is continuous on $(-1,0)$%
. Besides, $\lim_{x\rightarrow 0-}F(x)=-g_{1}(0)=0.$ This shows that $F$ is
continuous on $(-1,1)$, since $F\equiv 0$ on $[0,1).$ Combining all the
above arguments we conclude that $F\in C(\mathbb{R)}$.

Since $F\in C(\mathbb{R)}$ and $\Delta _{h_{1}...h_{k}}F_{2}(x)$ is a
continuous function of $x$, we obtain from (3.14) that $\Delta
_{h_{1}...h_{k}}F^{\ast }(x)$ is also a continuous function of $x.$ By Lemma
3.3, there exist a function $g_{2}\in C(\mathbb{R)}$ and a polynomial
function $H_{3}$ of order $k$ such that

\begin{equation*}
F^{\ast }=g_{2}+H_{3}.\eqno(3.17)
\end{equation*}%
It follows from (3.12), (3.14) and (3.17) that

\begin{equation*}
f=F+g_{2}+H_{2}+H_{3}.\eqno(3.18)
\end{equation*}

Introduce the notation

\begin{eqnarray*}
H(x) &=&H_{2}(x)+H_{3}(x)-H_{3}(0), \\
g(x) &=&F(x)+g_{2}(x)+H_{3}(0).
\end{eqnarray*}%
Obviously, $g\in C(\mathbb{R)}$ and $H(0)=0$. It follows from (3.18) and the
above notation that

\begin{equation*}
f=g+H.
\end{equation*}%
This completes the proof of the theorem.

\bigskip

\bigskip

\begin{center}
{\large \textbf{4. Ridge function representation}}
\end{center}

We start this section with the following lemma.

\bigskip

\textbf{Lemma 4.1.} \textit{Assume we are given pairwise linearly
independent vectors $\mathbf{a}^{i},$ $i=1,...,k,$ and a function $f\in C(%
\mathbb{R}^{n}) $ of the form}

\begin{equation*}
f(\mathbf{x})=\sum_{i=1}^{k}f_{i}(\mathbf{a}^{i}\cdot \mathbf{x}),\eqno(4.1)
\end{equation*}%
\textit{where $f_{i}$ are arbitrarily behaved univariate functions. Then for
any $h_{1},...,h_{k-1}\in \mathbb{R}$, and all indices $i=1,...,k$, $\Delta
_{h_{1}...h_{k-1}}f_{i}\in C(\mathbb{R})$.}

\bigskip

\textbf{Proof.} We prove this lemma for the function $f_{k}.$ It can be
proven for the other functions $f_{i}$ in the same way. Let $%
h_{1},...,h_{k-1}\in \mathbb{R}$ be given. Since the vectors $\mathbf{a}^{i}$
are pairwise linearly independent, for each $j=1,...,k-1,$ there is a vector
$\mathbf{b}^{j}$ such that $\mathbf{b}^{j}\cdot \mathbf{a}^{j}=0$ and $%
\mathbf{b}^{j}\cdot \mathbf{a}^{k}\neq 0$. It is not difficult to see that
for any $\lambda \in \mathbb{R}$, $\Delta _{\lambda \mathbf{b}^{j}}f_{j}(%
\mathbf{a}^{j}\cdot \mathbf{x})=0.$ Therefore, for any $\lambda
_{1},...,\lambda _{k-1}\in \mathbb{R}$, we obtain from (4.1) that

\begin{equation*}
\Delta _{\lambda _{1}\mathbf{b}^{1}...\lambda _{k-1}\mathbf{b}^{k-1}}f(%
\mathbf{x})=\Delta _{\lambda _{1}\mathbf{b}^{1}...\lambda _{k-1}\mathbf{b}%
^{k-1}}f_{k}(\mathbf{a}^{k}\cdot \mathbf{x}).\eqno(4.2)
\end{equation*}%
Note that in multivariate setting the difference operator $\Delta _{\mathbf{h%
}^{1}...\mathbf{h}^{k}}f(\mathbf{x})$ is defined similarly as in Section 2.
If in (4.2) we take

\begin{eqnarray*}
\mathbf{x} &\mathbf{=}&\frac{\mathbf{a}^{k}}{\left\Vert \mathbf{a}%
^{k}\right\Vert ^{2}}t\text{, }t\in \mathbb{R}\text{,} \\
\lambda _{j} &=&\frac{h_{j}}{\mathbf{a}^{k}\cdot \mathbf{b}^{j}}\text{, }%
j=1,...,k-1\text{,}
\end{eqnarray*}%
we will obtain that $\Delta _{h_{1}...h_{k-1}}f_{k}\in C(\mathbb{R})$.

\bigskip

Our main result is the following theorem.

\bigskip

\textbf{Theorem 4.1.} \textit{Assume a function $f\in C(\mathbb{R}^{n})$ is
of the form (4.1). Then there exist continuous functions $g_{i}:\mathbb{%
R\rightarrow R}$, $i=1,...,k$, and a polynomial $P(\mathbf{x})$ of degree at
most $k-1$ such that}

\begin{equation*}
f(\mathbf{x})=\sum_{i=1}^{k}g_{i}(\mathbf{a}^{i}\cdot \mathbf{x})+P(\mathbf{x%
}).\eqno(4.3)
\end{equation*}

\bigskip

\textbf{Proof.} By Lemma 4.1 and Theorem 3.1, for each $i=1,...,k$, there
exists a function $g_{i}\in C(\mathbb{R})$ and a polynomial function $H_{i}$
of $(k-1)$-th order with the property $H_{i}(0)=0$ such that

\begin{equation*}
f_{i}=g_{i}+H_{i}.\eqno(4.4)
\end{equation*}

Consider the function

\begin{equation*}
F(\mathbf{x})=f(\mathbf{x})-\sum_{i=1}^{k}g_{i}(\mathbf{a}^{i}\cdot \mathbf{x%
}).\eqno(4.5)
\end{equation*}%
It follows from (4.1), (4.4) and (4.5) that

\begin{equation*}
F(\mathbf{x})=\sum_{i=1}^{k}H_{i}(\mathbf{a}^{i}\cdot \mathbf{x}).\eqno(4.6)
\end{equation*}

Denote the restrictions of the multivariate functions $H_{i}(\mathbf{a}%
^{i}\cdot \mathbf{x})$ to the space $\mathbb{Q}^{n}$ by $P_{i}(\mathbf{x})$,
respectively. By Lemma 3.1, the functions $P_{i}(\mathbf{x})$ are ordinary
polynomials of degree at most $k-1$. Since the space $\mathbb{Q}^{n}$ is
dense in $\mathbb{R}^{n}$, and the functions $F(\mathbf{x})$, $P_{i}(\mathbf{%
x})$, $i=1,...,k$, are continuous on $\mathbb{R}^{n}$, and the equality

\begin{equation*}
F(\mathbf{x})=\sum_{i=1}^{k}P_{i}(\mathbf{x}),\eqno(4.7)
\end{equation*}%
holds for all $\mathbf{x}\in \mathbb{Q}^{n}$, we obtain that (4.7) holds
also for all $\mathbf{x}\in \mathbb{R}^{n}$. Now (4.3) follows from (4.5)
and (4.7) by putting $P=\sum_{i=1}^{k}P_{i}$.

\bigskip

Now we generalize Theorem 4.1 from $C(\mathbb{R}^{n})$ to any space $C^{s}(%
\mathbb{R}^{n})$ of $s$-th order continuously differentiable functions.

\bigskip

\textbf{Theorem 4.2.} \textit{Assume $f\in C^{s}(\mathbb{R}^{n})$ is of the
form (4.1). Then there exist functions $g_{i}\in C^{s}(\mathbb{R})$, $%
i=1,...,k$, and a polynomial $P(\mathbf{x})$ of degree at most $k-1$ such
that (4.3) holds.}

\bigskip

The proof is based on Theorem 4.1 and the following result of A. Pinkus \cite%
{Pinkus}.

\bigskip

\textbf{Theorem 4.3} (Pinkus \cite{Pinkus}). \textit{Assume $f\in C^{s}(%
\mathbb{R}^{n})$ is of the form (4.1). Assume, in addition, that each $%
f_{i}\in \mathcal{B}$. Then necessarily $f_{i}\in C^{s}(\mathbb{R)}$ for $%
i=1,...,k.$}

\bigskip

In Theorem 4.3, $\mathcal{B}$ denotes any linear space of real-valued
functions $u$ defined on $\mathbb{R}$, closed under translation, such that
if there is a function $v\in C(\mathbb{R)}$ for which $u-v$ satisfies the
Cauchy Functional Equation, then $u-v$ is necessarily linear, i.e. $%
u(x)-v(x)=cx,$ for some constant $c\in \mathbb{R}$.

Now the proof of Theorem 4.2 becomes obvious. Indeed, on the one hand, it
follows from Theorem 4.1 that the $s$-th order continuously differentiable
function $f-P$ can be expressed as $\sum_{i=1}^{k}g_{i}$ with continuous $%
g_{i}$. On the other hand, since the class $\mathcal{B}$ in Theorem 4.3, in
particular, can be taken as $C(\mathbb{R}),$ it follows that $g_{i}\in C^{s}(%
\mathbb{R})$.

\bigskip

\bigskip

\textbf{Remark 1.} Theorem 4.2 solves the problem posed in Buhmann and
Pinkus \cite{12} and Pinkus \cite[p. 14]{Pinkus1} up to a polynomial. In the
two dimensional setting $n=2$ it solves the problem completely. Indeed, it
is known that a bivariate polynomial $P(x,y)$ of degree $k-1$ is decomposed
into a sum of ridge polynomials with any given $k$ pairwise linearly
independent directions $(a_{i},b_{i}),$ $i=1,...,k$ (see e.g. \cite{97}).
That is,

\begin{equation*}
P(x,y)=\sum_{i=1}^{k}p_{i}(a_{i}x+b_{i}y),
\end{equation*}%
where $p_{i}$ are univariate polynomials of degree at most $k-1$.
Considering this in (4.3) gives the desired result.

\bigskip

\textbf{Remark 2.} Using our previous result \cite[Theorem 3.1]{A1}, the
degree of polynomial $P(\mathbf{x})$ in (4.3) can be reduced. Indeed, it
follows from (4.6) and (4.7) that the the above polynomial $P(\mathbf{x})$
is of the form (4.1). On the other hand, \cite[Theorem 3.1]{A1} states that
if a function $g\in C^{s}(\mathbb{R}^{n})$ is of the form (4.1), where $%
s\geq k-p+1$ and $p$ is the number of vectors $\mathbf{a}^{i}$ forming a
maximal linearly independent system, then there exist functions $g_{i}^{\ast
}\in C^{s}(\mathbb{R})$, $i=1,...,k$, and a polynomial $G(\mathbf{x})$ of
degree at most $k-p+1$ such that

\begin{equation*}
g(\mathbf{x})=\sum_{i=1}^{k}g_{i}^{\ast }(\mathbf{a}^{i}\cdot \mathbf{x})+G(%
\mathbf{x}).\eqno(4.8)
\end{equation*}%
Now putting $g(\mathbf{x})=P(\mathbf{x})$ in (4.8) and considering this in
(4.3) we see that our assertion is true.

\bigskip

\textbf{Remark 3.} In addition to the above $C^{s}(\mathbb{R})$, Theorem 4.1
can be restated also for the classes $C^{\infty }(\mathbb{R})$ of infinitely
differentiable functions and $D(\mathbb{R})$ of analytic functions. That is,
if under the conditions of Theorem 4.1, we have $f\in C^{\infty }(\mathbb{R}%
^{n})$ (or $f\in D(\mathbb{R}^{n})$), then this function can be represented
also in the form (4.3) with $g_{i}\in C^{\infty }(\mathbb{R})$ (or $g_{i}\in
D(\mathbb{R})$). This follows, similarly to the case $C^{s}(\mathbb{R})$
above, from Theorem 4.1 and Remark 2.2 in the book by Pinkus \cite{Pinkus1}.
In that remark, it was shown that, Theorem 4.3 can be restated for several
classes of functions, in particular, for the classes $C^{\infty }(\mathbb{R}%
) $ and $D(\mathbb{R})$.

\bigskip

The following corollaries show that for many directions $\mathbf{a}^{i}$, in
particular for those with rational coordinates, the polynomial terms in
Theorems 4.1 and 4.2 do not appear.

\bigskip

\textbf{Corollary 4.1.} \textit{Assume a function $f\in C(\mathbb{R}^{n})$
is of the form (4.1) and there is a nonsingular linear transformation $T:$ $%
\mathbb{R}^{n}\rightarrow \mathbb{R}^{n}$ such that $T\mathbf{a}^{i}\in
\mathbb{Q}\mathit{^{n}},$ $i=1,...,k$. Then there exist continuous functions
$g_{i}:\mathbb{R\rightarrow R}$, $i=1,...,k$, such that}

\begin{equation*}
f(\mathbf{x})=\sum_{i=1}^{k}g_{i}(\mathbf{a}^{i}\cdot \mathbf{x}).\eqno(4.9)
\end{equation*}

\bigskip

\textbf{Proof.} Applying the coordinate change $\mathbf{x\rightarrow y}$,
given by the formula $\mathbf{x}=T\mathbf{y}$, to both sides of (4.1) we
obtain that

\begin{equation*}
\tilde{f}(\mathbf{y})=\sum_{i=1}^{k}f_{i}(\mathbf{b}^{i}\cdot \mathbf{y}),
\end{equation*}%
where $\tilde{f}(\mathbf{y})=f(T\mathbf{y})$ and $\mathbf{b}^{i}=T\mathbf{a}%
^{i},$ $i=1,...,k.$ Let us repeat the proof of Theorem 4.1 for the function $%
\tilde{f}$. Since the vectors $\mathbf{b}^{i}$, $i=1,...,k,$ have rational
coordinates, it is not difficult to see that the restrictions of the
functions $H_{i}$ to $\mathbb{Q}$ are univariate polynomials. Indeed, for
each $\mathbf{b}^{i}$ we can choose a vector $\mathbf{c}^{i}$ with rational
coordinates such that $\mathbf{b}^{i}\cdot \mathbf{c}^{i}=1$. If in the
equality $H_{i}(\mathbf{b}^{i}\cdot \mathbf{x})=P_{i}(\mathbf{x}),$ $\mathbf{%
x}\in \mathbb{Q}^{n}$, we take $\mathbf{x=c}^{i}t$ with $t\in \mathbb{Q}$,
we obtain that $H_{i}(t)=P_{i}(\mathbf{c}^{i}t)$ for all $t\in \mathbb{Q}$.
Now since $P_{i}$ is a multivariate polynomial on $\mathbb{Q}^{n}$, $H_{i}$
is a univariate polynomial on $\mathbb{Q}$. Denote this univariate
polynomial by $L_{i}$. Thus the formula

\begin{equation*}
P_{i}(\mathbf{x})=L_{i}(\mathbf{b}^{i}\cdot \mathbf{x})\eqno(4.10)
\end{equation*}%
holds for each $i=1,...,k$, and all $\mathbf{x}\in \mathbb{Q}^{n}$. Since $%
\mathbb{Q}^{n}$ is dense in $\mathbb{R}^{n}$, we see that (4.10) holds, in
fact, for all $\mathbf{x}\in \mathbb{R}^{n}$. Thus the polynomial $P(\mathbf{%
x})$ in (4.3) can be expressed as $\sum_{i=1}^{k}L_{i}(\mathbf{b}^{i}\cdot
\mathbf{x})$. Considering this in Theorem 4.1, we obtain that

\begin{equation*}
\tilde{f}(\mathbf{y})=\sum_{i=1}^{k}g_{i}(\mathbf{b}^{i}\cdot \mathbf{y}),%
\eqno(4.11)
\end{equation*}%
where $g_{i}$ are continuous functions. Using the inverse transformation $%
\mathbf{y}=T^{-1}\mathbf{x}$ in (4.11) we arrive at (4.9).

\bigskip

\textbf{Corollary 4.2.} \textit{Assume a function $f\in C^{s}(\mathbb{R}^{n})
$ is of the form (4.1) and there is a nonsingular linear transformation $T:$
$\mathbb{R}^{n}\rightarrow \mathbb{R}^{n}$ such that $T\mathbf{a}^{i}\in
\mathbb{Q}\mathit{^{n}},$ $i=1,...,k$. Then there exist functions $g_{i}\in
C^{s}(\mathbb{R})$, $i=1,...,k$, such that}

\begin{equation*}
f(\mathbf{x})=\sum_{i=1}^{k}g_{i}(\mathbf{a}^{i}\cdot \mathbf{x}).
\end{equation*}

\bigskip

The proof of Corollary 4.2 easily follows from Corollary 4.1 and Theorem 4.3.

\bigskip

\textbf{Acknowledgments.} The original version of this paper had a gap in
the proof of the main result. Namely, Lemma 3.1 was incorrectly applied. 
This was found by Professor A.A. Kuleshov and
the anonymous reviewer, to whom the authors are very grateful.

\end{document}